\documentclass[12pt]{article}
\usepackage{amsfonts,enumerate,amsmath,amssymb}

\newtheorem{theorem}{Theorem}

\def\CC{{\mathbb C}}
\def\RR{{\mathbb R}}
\def\ZZ{{\mathbb Z}}
\def\QQ{{\mathbb Q}}

\begin{document}

\title{The generalized Kummer construction and the cohomology rings of $G_2$-manifolds
\thanks{Th work is supported by RSCF (grant No 14-11-00441).}
}
\author{I.A. Taimanov
\thanks{Sobolev Institute of Mathematics, 630090 Novosibirsk, and Novosibirsk State University, 630090 Novosibirsk, Russia;
e-mail: taimanov@math.nsc.ru}
}
\date{}
\maketitle

First examples of closed Riemannian manifolds with holonomy groups
$G_2$ and $Spin(7)$ were constructed by Joyce by using the generalized Kummer construction
\cite{Joyce1,Joyce2,Joyce3} (see also \cite{Joyce4}),
which in its classical version is applied to the $4$-torus for constructing
$K3$ surfaces. These manifolds are obtained by resolution of singularities
of the orbifolds $T^n/\Gamma$, where $\Gamma$ acts discretely on the torus
 ($n=4, \Gamma = \ZZ_2$ for $K3$ surfaces, $n=7$ for $G_2$-manifolds and $n=8$ for $Spin(7)$-manifolds). In \cite{Joyce5} this construction was generalized for the orbifolds $M^8/\Gamma$, where $M^8$ is an $8$-dimensional
 Calabi--Yau manifold, for constructing new examples of $Spin(7)$ manifolds. The manifolds obtained in these articles together with $G_2$-manifolds constructed by using twisted connected sums \cite{Kovalev,KL,CHNP} give all known examples of closed
$G_2$- and $Spin(7)$-manifolds.

In \cite{CHNP} the cohomology rings were calculated for many $G_2$-manifolds obtained as twisted connected sums and
it was noted that the cohomology rings for $G_2$-manifolds obtained by using the generalized Kummer construction
are still not calculated.

The ring structure in the cohomology of  a $K3$ surface was established by Milnor
by using rather general algebraic arguments \cite{Milnor} and the explicit form of the basis of two-cycles
for which the intersection form is canonical was unknown until recently.
In \cite{Taimanov} we found this basis by using the intersection theory on manifolds.

In this article we demonstrate how the intersection theory can be applied to calculation of the cohomology rings of
$G_2$-manifolds constructed in \cite{Joyce1,Joyce2}.
Since these manifolds were obtained as a finite family of examples we shall make calculations for a typical example, for other cases the cohomology rings are calculated analogously. Therewith we restrict ourselves to the rational cohomology rings.
A detailed investigation of the cohomology rings with integer coefficients can be done as in
\cite{Taimanov} but will involve bulky calculations.

We present the results in terms of the intersection ring $H_\ast(X;\QQ)$ which is dual to the cohomology ring
$H^\ast(X;\QQ)$.

For a closed oriented $n$-dimensional smooth manifold $X$, an intersection of cycles
$$
H_k (X;\ZZ) \times H_l(X;\ZZ) \stackrel{\cap}{\longrightarrow} H_{k+l-n}(X;\ZZ)
$$
is defined as follows
$$
u \cap v = D^{-1} (Du \cup Dv),
$$
where
$$
D: H_i(X;\ZZ) \to H^{n-i}(X;\ZZ), \ \ \ i=0,\dots,n,
$$
is the Poincare duality. If the cycles $u \in H_k(X;\ZZ)$ and $v \in H_l(X;\ZZ)$ are realized by transversally intersecting smooth submanifolds $Y$ and $Z$, then their intersection is the smooth $(k+l-n)$-dimensional submanifold $W$ which realizes the cycle
$w \in H_{k+l-n}(X;\ZZ)$ such that
$$
u \cap v = w.
$$
The orientation of the intersection $Y \cap Z$ is defined as follows.
Let $x \in Y \cap Z$, $(e_1,\dots,e_{k+l-n})$ be a basis in the tangent space to  $Y \cap Z$ at $x$ and
$(e_1,\dots,e_{k+l-n},e^\prime_1,\dots,e^\prime_{n-l})$ and  $(e_1,\dots,e_{k+l-n},e^{\prime\prime}_1,\dots,
e^{\prime\prime}_{n-k})$ be positively ori\-ented bases in the tangent spaces to $Y$ and $Z$ at $x$. Then $Y \cap Z$ is oriented such that the basis $(e_1,\dots,e_{k+l-n})$
in the tangent space to $Y \cap Z$ determines the same orientation (positive or negative) as the basis
$(e_1,\dots,e_{k+l-n}$, $e^\prime_1,\dots$, $e^\prime_{n-l}$, $e^{\prime\prime}_1$, $\dots$,
$e^{\prime\prime}_{n-k})$ in the tangent space $X$. Therefore the intersection of homologies satisfies
the natural anticommutativity condition
$$
u \cap v = (-1)^{(n-k)(n-l)} v \cap u, \ \ \ \ \ \dim u = k, \dim v = l, \dim M = n.
$$

Since not all cycles are realized by smooth submanifolds, in general the geometric definition of intersection needs an introduction of a special class of chains and that was done by Lefschetz who generalized the Poincare theory of intersection of cycles of complemented $k$ and $n-k$ dimensions onto cycles of any dimensions. A strong and complete exposition of this construction was done in \cite{GP} where the homology group $H_\ast(X)$ with the intersection operation is called the Lefschetz ring
and its topological invariance is established by using its isomorphism $D$ to $H^\ast(X)$.

Let us turn to a concrete example of a Joyce manifold.

Let $T^7 = \RR^7/\ZZ^7$ be the $7$-torus which is the quotient space of $\RR^7$
under the action of the integer lattice $\ZZ^7$ by translations $x \to x + y$ where $x \in \RR^7, y \in \ZZ^7 \subset \RR^7$.
We denote by $x_1,\dots,x_7$ the coordinates, on $T^7$, induced by the Euclidean coordinates on $\RR^7$ and defined up to integers.
Let
\begin{equation}
\label{involutions}
\begin{split}
\alpha((x_1,\dots,x_7)) = (-x_1,-x_2,-x_3,-x_4,x_5,x_6,x_7),
\\
\beta((x_1,\dots,x_7)) = (b_1-x_1,b_2-x_2,x_3,x_4,-x_5,-x_6,x_7),
\\
\gamma((x_1,\dots,x_7)) = (c_1-x_1,x_2,c_3-x_3,x_4,c_5-x_5,x_6,-x_7),
\end{split}
\end{equation}
be transformations, of $T^7$, which are parameterized by constants $b_1$, $b_2$, $c_1$, $c_3$, $c_5$.
They are involutions
$$
\alpha^2 = \beta^2 = \gamma^2 = 1,
$$
and pairwise commute:
$$
\alpha\beta=\beta\alpha, \ \ \alpha\gamma = \gamma\alpha, \ \ \ \beta\gamma = \gamma\beta,
$$
and therefore for all $b_1,b_2,c_1,c_3,c_5$ generate an action of $\Gamma=\ZZ_2^3$ on $T^7$.

We consider the following constants which lead to a simply-connected manifold $M^7$
(Example 3 in \cite{Joyce2}):
\begin{equation}
\label{data1}
b_1 = c_5 = 0, \ \ \ b_2 = c_1 = c_3 = \frac{1}{2}.
\end{equation}

Let us recall the main facts on the topology of the $G_2$-manifold $M^7$ which is constructed from these data
(\ref{involutions}) and (\ref{data1})
\cite{Joyce1,Joyce2}:

\begin{enumerate}
\item
$\Gamma$ acts on $H^\ast(T^7)$ by involutions, there no nontrivial invariant subspaces in
$H^1(T^7)$ and $H^2(T^7)$, and the invariant subspace in $H^3(T^7)$ is of dimension $7$ and is generated by the cohomology classes induced by the following forms with integer periods:
\begin{equation}
\label{invforms}
\begin{split}
dx_2 \wedge dx_4 \wedge dx_6, \ \ \
dx_3 \wedge dx_4 \wedge dx_7, \ \ \
dx_5 \wedge dx_6 \wedge dx_7,
\\
dx_1 \wedge dx_2 \wedge dx_7, \ \ \
dx_1 \wedge dx_3 \wedge dx_6,
\\
dx_1 \wedge dx_4 \wedge dx_5, \ \ \
dx_2 \wedge dx_3 \wedge dx_5.
\end{split}
\end{equation}
This implies that
$$
b^1(T^7/\Gamma) = b^2(T^7/\Gamma) = 0, \ \
b^3(T^7/\Gamma) =  7,
$$
where by
$$
b^k(Y) = \dim_\QQ H^k(Y;\QQ), \ \ \ k\geq 0,
$$
we denote the Betti numbers of $Y$.
Since the $7$-form $dx_1 \wedge \dots \wedge dx_7$ is $\Gamma$-invariant,
the Hodge operator $\ast: H^k(T^7) \to H^{7-k}(T^7)$ maps invariant forms into invariant forms and therefore
$$
b^6(T^7/\Gamma) = b^5(T^7/\Gamma) = 0, \ \ \ b^4(T^7/\Gamma) = 7.
$$

\item
The action of $\Gamma = \ZZ_2^3$ is not free. The fixed points of every involution form the family $\alpha, \beta, \gamma$ form $16$ fixed $3$-tori:
$$
(\vartheta^\alpha_1,\vartheta^\alpha_2,\vartheta^\alpha_3,\vartheta^\alpha_4, x_5,x_6,x_7),  \ \ \
 \vartheta^\alpha_1,\dots,\vartheta^\alpha_4 \in
\left\{0, \frac{1}{2}\right\} \ \ \ \ \mbox{for $\alpha$};
$$
$$
(\vartheta^\beta_1,\vartheta^\beta_2,x_3,x_4,\vartheta^\beta_5,\vartheta^\beta_6,x_7), \ \ \
\ \vartheta^\beta_1, \vartheta^\beta_5, \vartheta^\beta_6 \in \left\{0, \frac{1}{2}\right\}, \
\vartheta^\beta_2 \in \left\{ \frac{1}{4} ,\frac{3}{4}\right\} \ \ \ \ \mbox{for $\beta$};
$$
$$
(\vartheta^\gamma_1,x_2, \vartheta^\gamma_3, x_4,\vartheta^\gamma_5, x_6, \vartheta^\gamma_7),
\ \vartheta^\gamma_5, \vartheta^\gamma_7 \in \left\{0, \frac{1}{2}\right\}, \  \ \
\vartheta^\gamma_1, \vartheta^\gamma_3 \in \left\{ \frac{1}{4}, \frac{3}{4}\right\}
\ \ \ \ \mbox{for $\gamma$}.
$$

\item
$\Gamma/\ZZ_2$ acts by nontrivial transpositions on the fixed tori of the involutions $\alpha, \beta, \gamma$ and
the $\Gamma$-orbit every torus fixed by certain involution consists of
$4$ tori. By parameterizing the tori by the quadruples of parameters $\vartheta_i$, we describe the $\Gamma$-orbits of
tori fixed by involutions $\alpha, \beta$, or $\gamma$:

the orbits of tori fixed by $\alpha$ splits into two pairs corresponding to different values of $\vartheta^\alpha_4 \in \left\{0,\frac{1}{2}\right\}$:
$$
\left\{(0,0,0,\vartheta^\alpha_4), \left(0,\frac{1}{2},0,\vartheta^\alpha_4\right), \left(\frac{1}{2},0,\frac{1}{2},\vartheta^\alpha_4\right),\left(\frac{1}{2},
\frac{1}{2},\frac{1}{2},\vartheta^\alpha_4\right)\right\},
$$
$$
\left\{\left(0,0,\frac{1}{2},\vartheta^\alpha_4\right),\left(0,\frac{1}{2},\frac{1}{2},
\vartheta^\alpha_4\right),\left(\frac{1}{2},0,0,\vartheta^\alpha_4\right),\left(\frac{1}{2},\frac{1}{2},0,\vartheta^\alpha_4\right)\right\};
$$

the orbits of tori fixed by $\beta$ are parameterized by $\vartheta^\beta_5, \vartheta^\beta_6 \in \left\{0, \frac{1}{2}\right\}$;

the orbits of tori fixed by $\gamma$ are parameterized by $\vartheta^\gamma_5, \vartheta^\gamma_5 \in \left\{0, \frac{1}{2}\right\}$.

\item
The products, of the elementary involutions, $\alpha\beta, \beta\gamma, \alpha\gamma$ and $\alpha\beta\gamma$
have no fixed points because they act on some of coordinates by translations: for instance,
$\alpha\beta(x_2) = x_2 + \frac{1}{2}$.

\item
Each of the involutions $\alpha,\beta,\gamma$ generates such a $\ZZ_2$-action on $T^7$ that
$$
T^7/\ZZ_2 = T^3 \times (T^4/\ZZ_2),
$$
where $T^4/\ZZ_2$ is the Kummer surface from which a $K3$ surface is obtained by resolving $16$ singular points (see, for instance,
\cite{AS,H}).
All generators of $\pi_1(T^7)$ under the action of certain involution are projected into a contractible contour in
$T^4/\ZZ_2$. Therefore $\pi_1(T^7/\Gamma) = 0$.

\item
The singularity set of the orbifold $T^7/\Gamma$ splits into twelve $3$-tori corresponding to
the $\Gamma$-orbits of fixed tori of the involutions $\alpha, \beta, \gamma$.
For every singular torus there is a neighborhood of which is homeo\-mor\-phic to
\begin{equation}
\label{singularity}
U = T^3 \times (D/\ZZ_2),
\end{equation}
where $D = \{|z|\leq \tau \ : \ z \in \CC^2\}$ is a small neighborhood of the origin in $\CC^2$ and
$\ZZ_2$ acts on $D$ as follows: $z\to -z$.

\item
From a topological point of view the manifold $M^7$
is constructed from $T^7/\Gamma$ by a fiber-wise resolution of singularities (\ref{singularity})
in $D/\ZZ_2$. For every singular torus we remove the neighborhood $U$ and replace it by the product
$$
\widetilde{U} = T^3 \times V,
$$
where $V$ is a sufficiently small neighborhood of the zero section
of the fiber bundle
$\gamma^2 \to \CC P^1$ and $\gamma \to \CC P^1$ is the tautological bundle over $\CC P^1$.
We have the fibration
\begin{equation}
\label{fibration}
\widetilde{\pi}: \widetilde{U} = T^3 \times V \to T^3.
\end{equation}
Therewith the zero section of the fiber bundle $U \to T^3$ is replaced by $\CC P^1$, and analogously to the case
of the resolution of singularities of the Kummer surface $T^4/\ZZ_2$ the selfintersection index of the
homology class $[\CC P^1]$ in the four-dimensional fiber is equal to $-2$. The homology group of $T^3 \times V$ is equal to
$$
H_\ast (T^3 \times V) = H_\ast (T^3) \otimes H_\ast (\CC P^1).
$$

\item
The rational homology groups of $M^7$ are as follows:

$b^2 = 12$ and the generators are given $12$ cycles which are determined by the submanifolds diffeomorphic to
$\CC P^1$ and corresponding to the resolutions of $12$ singular tori;

$b^3 = 43$ and the generators are presented by $7$ cycles realized by the tori corresponding to
the invariant $3$-forms on $T^7$ and by $12$ families of the products of $\CC P^1$ and the generators of $1$-cycles
on the singular tori.

The Betti numbers $b^4 = 43$ and $b^5=12$ are found via the Poincare duality.
\end{enumerate}

Let us construct the submanifolds realizing the generators of $H_\ast(M^7;\QQ)$.

{\sc $t$-cycles.}
Let us take three-tori $T_\alpha, T_\beta, T_\gamma$ such that they are homologous to the fixed tori of the involutions
$\alpha, \beta, \gamma$ but their $\Gamma$-orbits consists of $8$ tori.
For instance, for $\alpha$ we consider the torus whose points have the following coordinates
$$
T_\alpha = \left\{ \left(\frac{1}{3},\frac{1}{3},\frac{1}{3},\frac{1}{3}, x_5,x_6,x_7\right), \ \ x_5,x_6,x_7\in \RR\right\}.
$$
The fixed tori which are homologous to it we enumerate them in some way and
denote by $T_{\alpha 1}, \dots, T_{\alpha 4}$.
Let us construct a $4$-torus  $T^\prime_\alpha$, dual to $T_\alpha$,
such that the intersection index of $T^\prime_\alpha$ and $T_\alpha$ is equal to
$\mathrm{ind} = 1$:
$$
T_\alpha \cap T^\prime_\alpha = \mathrm{pt},
$$
and its $\Gamma$-orbit consists of $8$ tori. For instance, we may put
$$
T^\prime_\alpha = \left\{x_1,x_2,x_3,x_4, \left(\frac{1}{3},\frac{1}{3},\frac{1}{3}\right), \ \ x_1,x_2,x_3,x_4 \in \RR\right\}.
$$
Analogously we construct such tori corresponding to the involutions $\beta$ and $\gamma$.
By construction, for every such a torus an integral of certain from from the first line
(\ref{invforms}) over the torus is equal to $1$, and the integrals of other forms from
(\ref{invforms}) over this torus vanish: for instance, $\int_{T_\alpha} dx_5 \wedge dx_6 \wedge dx_7 =1$.
For the remaining forms from (\ref{invforms}) we construct the analogous
$3$-tori $T_1,\dots,T_4$ and the dual $4$-tori $T^\prime_1,\dots,T^\prime_4$.

The $\Gamma$-orbits of the dual tori  $T_l$ and $T^\prime_l$ intersect at $64 = 8^2$ points and
the projections of these tori into $T^7/\Gamma$ intersect at $8$ points which do not lie in the singular tori.
By the intersection theory in $T^7$, for $k\neq l$ the tori $T_k$ and $T_l^\prime$ intersect in the cycle homologous to zero, where $k,l \in \{\alpha,\beta,\gamma,1,2,3,4\}$.

Let us consider the projections of $T_k$ and $T^\prime_l$ onto $T^7/\Gamma$. Denote the cycles which represent these projections by $t_k=[\pi(T_k)]$ and $t^\prime_l=[\pi(T^\prime_l)]$, where $\pi: T^7 \to T^7/\Gamma$ is the projection onto the quotient.
We have
\begin{equation}
\label{rel1}
t_k \cap t^\prime_l =
\begin{cases} 8 & \mbox{for $k=l$, where $k,l \in \{\alpha,\beta,\gamma,1,2,3,4\}$}\\
0 & \mbox{otherwise}.
\end{cases}
\end{equation}

The projections of tori which represent the cycles $t_k$ and $t^\prime_l$ where $k \in \{\alpha,\beta,\gamma,1,2,3,4\}$ and $l \in \{1,2,3,4\}$ do not intersect the singular tori and therefore they are not affected by the resolution of singularities. We denote the corresponding cycles in $M^7$ by the same symbols.

The tori $T^\prime_\alpha,T^\prime_\beta,T^\prime_\gamma$ intersect the fixed tori and after the projection into  $T^7/\Gamma$ and the resolution of singularities they are transformed into the cycles realized by
$K3$ surfaces. For instance, $\alpha$ acts on $T^\alpha$ as a reflection with $16$ fixed points
at which $T^\alpha$ intersects tori fixed by $\alpha$. The projection
$T^7 \to T^7/\langle \alpha\rangle$, where $\langle\alpha\rangle = \ZZ_2$ is generated by
$\alpha$, maps $T^\prime_\alpha$ into a Kummer surface on which $\ZZ_2^2 = \Gamma/\langle \alpha \rangle$
acts by shifts.  Finally $T^\prime_\alpha$ is mapped by the projection
$T^7 \to T^7/\Gamma$ into a Kummer surface which intersect every singular torus at four points.
After the resolution of singularities the Kummer surface goes to a $K3$ surface. Analogously we transform
$T^\prime_\beta$ and $T^\prime_\gamma$.
The cycles, in $H_4(M^7)$, realized by these submanifolds we denote by $t^\prime_\alpha,t^\prime_\beta, t^\prime_\gamma$.

{\sc $C$-cycles.}
For every involution
$\delta \in \{\alpha,\beta,\gamma\}$ we enumerate the singular tori $T_{\delta i}$ by numbers from $1$ to $4$.
We denote by
$$
C_{\alpha 1}, \dots, C_{\alpha 4}, C_{\beta 1}, \dots, C_{\gamma 4} \ \ \mbox{($\dim =2$)}
$$
the two-manifolds which are diffeomorphic to $\CC P^1$, lie in the fibers of fibrations (\ref{fibration})
and generate the nontrivial homology of $V$.

For every singular torus $T_{\delta i}$ we take three nontrivial loops $\lambda_{\delta ij}$, $\delta \in \{\alpha, \beta,\gamma\}$, $i=1,\dots,4, j=1,2,3$, which lie on the torus, start at a certain fixed point and are drawn via varying one of the variables $x_i$.  For every such a loop we choose a $2$-torus $\tau_{\delta ij} \subset T_{\delta i}$ for which the intersection index with the loop is equal to $\mathrm{ind} =1$. Here we assume that the orientation on $T_{\delta i}$
is determined by the corresponding from from the first line in (\ref{invforms}). Let us construct the following submanifolds:
$$
\lambda_{\delta i j} \times C_{\delta i} \  \ \mbox{($\dim =3$)}, \ \ \
\tau_{\delta ij} \times C_{\delta i} \ \ \mbox{($\dim=4$)}, \ \
T_{\delta i} \times C_{\delta i} \ \ \mbox{($\dim = 5$)}.
$$
Here we mean by the product of $W$, a subset of the torus,
and $C_{\delta i}$ the product $W$ and $C_{\delta i} \subset V$.
Let us denote the cycles realized by these submanifolds by
$$
c_{\delta i} = [C_{\delta i}], \ \ \ c_{\delta ij} = [\lambda_{\delta i j} \times C_{\delta i}], \ \
c_{\delta ij}^\prime = [\tau_{\delta ij} \times C_{\delta i}], \ \  c_{\delta i}^\prime = [T_{\delta i} \times C_{\delta i}].
$$
We find the intersections of these cycles.

a) It is evident that the intersection of $c$-cycles with different first indices vanish and also
the intersection of every $c$-cycle with $t_k$, $k \in \{\alpha,\beta,\gamma,1,2,3,4\}$, $t^\prime_l$, $l =1,2,3,4$,
or $t_\delta$, $\delta \in \{\alpha,\beta,\gamma\}$, is equal to zero.

b) If the sum of the dimensions of cycles is less than the dimension of the manifold, then by small perturbation the submanifolds realizing the cycles are transformed into non-intersecting submanifolds. Therefore the product of such cycles vanishes:
$$
c_{\delta i} \cap c_{\delta i} =0, \ \ c_{\delta i} \cap c_{\delta ij} = 0, \ \ c_{\delta i} \cap c^\prime_{\delta ij} = 0, \ \
c_{\delta ij} \cap c_{\delta ik} =0,
$$

c) In a $3$-torus the tori $\tau_{\delta ij}$ may be shifted into parallel tori and that implies
$$
c^\prime_{\delta ij} \cap c^\prime_{\delta ij} = 0.
$$
Also, if the loop $\lambda_{\delta ij}$ lies in $\tau_{\delta ijk}$, $i \neq k$, then
it can be shifted from the torus with eliminating an intersection and that implies
$$
c_{\delta ij} \cap c_{\delta ik}^\prime  = 0 \ \ \ \mbox{for $j \neq k$}.
$$

d) Low-dimensional tori (loops and $2$-tori) in the fixed tori are homologous to low-dimensional subtori, in  $T_k$,
which are homologous to zero over the rational field of coefficients, because on $T^7$ there are no $\Gamma$-invariant $1$-
and $2$-forms. Let us consider the intersections of the cycles $ [\tau_{\delta ij}\times C_{\delta i}]$ and
$ [T_{\delta i} \times C_{\delta i}]$. By perturbing submanifolds homologous to $C_{\delta i}$ in every fiber,
we obtain the intersection set of the form
$\pm 2 \tau_{\delta ij}$,
which is homologous to zero. Therefore
$$
c_{\delta ij}^\prime \cap  c_{\delta i}^\prime = 0.
$$
Analogously it is proved that
$$
c_{\delta ij}^\prime \cap c_{\delta ik}^\prime = 0 \ \ \mbox{for $j\neq k$}, \ \
c_{\delta ij} \cap c^\prime_{\delta i} = 0,
$$
because in this case the intersection set is homologous to $\pm 2  \lambda$, where $\lambda$ is a loop in
$T_{\delta i}$.

e) Since in $T_{\delta i}$ $\lambda_{\delta ij}$ and
$\tau_{\delta ij}$ intersect at a point $P$ with index equal to one, by small perturbation
the intersection of $\lambda_{\delta ij}\times C_{\delta i}$ and $\tau_{\delta ij} \times C_{\delta i}$
is reduced to the selfintersection of $\CC P^1$ in the fiber of (\ref{fibration}) over $P$ with the selfintersection index equal to
$\mathrm{ind} =-2$.
Therefore
$$
c_{\delta ij} \cap c^\prime_{\delta ij} = -2.
$$

f) Analogously to the previous case the selfintersection of $c^\prime_{\delta i}$ reduces to
the selfintersections $\CC P^1$ in fibers on (\ref{fibration}) over all points of $T_{\delta i}$ and therefore after small perturbation
of two copies of
$T_{\delta i} \times C_{\delta i}$ we obtain a transversal intersection homologous to
$-2 T_{\delta i}$.
Therefore
$$
c^\prime_{\delta i} \cap c^\prime_{\delta i} = -2 t_\delta.
$$

We calculated the ring structure in $H_\ast(M^7;\QQ)$:

\begin{theorem}
The rational homology group $H_\ast(M^7;\QQ)$ of the $G_2$-manifold obtained by using
the generalized Kummer construction from the data (\ref{involutions}) and (\ref{data1}),
has the following generators of dimension $\leq \dim M^7=7$:
$$
\dim =2: \ \ c_{\delta i} ; \ \ \ \dim = 3:  \ \  c_{\delta ij}, \ t_\delta, \ t_i;
$$
$$
\dim = 4: \ \  c_{\delta ij}^\prime, \ t_\delta^\prime, \ t_i^\prime; \ \ \ \
\dim = 5: \ \ c_{\delta i}^\prime,
$$
where $\delta \in\{\alpha,\beta,\gamma\}, i=1,\dots,4, j=1,2,3$.
All these cycles are realized by embedded submanifolds and hence by integral cycles, i.e., lie in
$H_\ast(M^7;\ZZ) \subset H_\ast(M^7;\QQ)$.

The intersection operation is given by nontrivial intersections which are (up to the anticommutativity)
\begin{equation}
\label{relations}
\begin{split}
c_{\delta i} \cap c_{\delta i}^\prime = -2, \ \  c_{\delta ij} \cap c_{\delta ij}^\prime = -2, \ \
t_\delta \cap t_\delta^\prime =8, \ \ t_i \cap t_i^\prime = 8,
\\
c_{\delta i}^\prime \cap c_{\delta i}^\prime = -2 t_\delta.
\end{split}
\end{equation}
\end{theorem}

{\sc Remarks.}
1) For describing the rational cohomology rings it is enough to
to correspond to any generator $a \in H_\ast, \dim a =k$
a generator $\bar{a}\in H^\ast, \deg \bar{a} = n-k$ and replace all relations of the form $a \cap b = c$ by relations
$\bar{a} \cup \bar{b} = \bar{c}$.

2) Analogously one can describe the rational cohomology rings of other manifolds, with special holonomy, obtained by using the generalized Kummer construction.

3) The first line in (\ref{relations}) means that there is a pairing between $H_k$ and $H_{n-k}$ (the Poincare duality) which in these generators has the form $a \to a^\prime$ and is defined by the diagonal matrix
$\det = 2^{12}$ for $k=2$ and $\det = 2^{36} 8^{7} = 2^{57}$ for $k=3$. These cycles do not generate
$H_\ast(M^7;\ZZ)/\mathrm{Torsion}$, because for the generators of  this group the pairings are defined by unimodular matrices.
Nontrivial products are given in the second line of (\ref{relations}).
One can obtain generators of $H_\ast(M^7;\ZZ)/\mathrm{Torsion}$ by explicit but probably bulky geometric constructions analogously to how that was done for $K3$ surfaces in
\cite{Taimanov}.

4) The Deligne--Griffiths--Morgan--Sullivan theorem \cite{DGMS} reads that simply-connected closed K\"ahler manifolds are formal \cite{Sullivan}.
This result is not gene\-ra\-lized for symplectic manifolds as it was shown first for manifolds of dimension $\geq 10$ \cite{BT1} (see also \cite{BT2}) and later for the remaining dimension $8$ \cite{FM} (simply-connected closed manifolds of
dimension $\leq 6$ are always formal). If we treat K\"ahler manifolds as manifolds with special holonomy $U(n)$, then
naturally appears the conjecture of formality of simply-connected closed manifolds with special holonomy. Manifolds with
$SU(n)$ and $Sp(n)$ holonomy are evidently formal because they are K\"ahler.
For the remaining cases of quaternion--K\"ahler manifolds (with holonomy $Sp(n)Sp(1)$) and with holonomy groups
$G_2$ (of dimension $7$) and $Spin(7)$ (of dimension $8$) counterexamples to this conjecture are unknown.

The only obstructions to the formality of seven-dimensional manifolds are nontrivial triple Massey products
of the form $\langle a,b,c \rangle$, where $a,b,c \in H^2$ and
$a \cup b = b \cup c=0$. (For $b_1=1$ there could be only one obstruction which is the triple Massey product $\langle x,x,x\rangle$ with $x$ a generator of $H^2$. However this Massey produc vanishes due to the reasons of dimension, and such a manifold is formal \cite{FIM}). Theorem 1 implies that there are many sets $a,b,c$ of two-dimensional cohomology classes which meet the latter conditions but we did not manage to construct from them a nontrivial Massey product.

\end{document}